\documentclass[sn-mathphys]{sn-jnl}

\jyear{2024}%

\theoremstyle{thmstyleone}%
\theoremstyle{thmstyletwo}%

\theoremstyle{thmstylethree}%

\usepackage{xcolor}

\newcommand\new[1]{\textcolor{black}{#1}}
\newcommand\added[1]{\new{#1}}
\usepackage{soul}
\newcommand\deleted[1]{}
\newcommand\replace[2]{\deleted{#1}\new{#2}}
\newcommand\reviewcomment[1]{}

\begin{document}

\title{Adaptive time step selection for Spectral Deferred Correction\deleted{s}}


\author*[1,2]{\fnm{Thomas} \sur{Saupe}}\email{t.baumann@fz-juelich.de}
\author[2]{\fnm{Sebastian} \sur{G\"otschel}}
\author[2]{\fnm{Thibaut} \sur{Lunet}}
\author[2]{\fnm{Daniel} \sur{Ruprecht}}
\author[1]{\fnm{Robert} \sur{Speck}}

\affil[1]{\orgdiv{J{\"u}lich Supercomputing Centre}, \orgname{Forschungszentrum J{\"u}lich GmbH}, \orgaddress{\street{Wilhelm-Johnen-Stra{\ss}e}, \city{J{\"u}lich}, \postcode{52428}, \country{Germany}}}

\affil[2]{\orgdiv{Chair Computational Mathematics}, \orgname{Institute of Mathematics, Hamburg University of Technology}, \orgaddress{\street{Am Schwarzenberg-Campus 3}, \city{Hamburg}, \postcode{21073}, \country{Germany}}}

\abstract{
Spectral Deferred Correction\deleted{s} (SDC) is an iterative method for the numerical solution of ordinary differential equations. It works by refining the numerical solution for an initial value problem by approximately solving differential equations for the error, and can be interpreted as a preconditioned fixed-point iteration for solving the fully implicit collocation problem.
We adopt techniques from embedded Runge-Kutta Methods (RKM) to SDC in order to provide a mechanism for adaptive time step size selection and thus increase computational efficiency of SDC.
We propose two SDC-specific estimates of the local error that are generic and \replace{require only minimal problem specific tuning}{do not rely on problem specific quantities}.
We demonstrate a gain in efficiency over standard SDC with fixed step size \replace{,}{and} compare efficiency favorably against state-of-the-art adaptive RKM\deleted{ and show that due to its iterative nature, adaptive SDC can cope efficiently with silent data corruption}.
}




\maketitle

\section{Introduction}\label{secIntroduction}
Spectral Deferred Corrections (SDC) were introduced by Dutt et al.~\cite{DuttSDCOriginal} as a more stable variant of the classical deferred correction\deleted{s} approach for solving ordinary differential equations (ODEs).
Deferred correction\deleted{s} methods iteratively refine the numerical solution for an initial value problem by approximately solving differential equations for the error.

While SDC often requires more work per time step than classical 
Runge-Kutta Methods (RKM)\footnote{\added{Note that SDC with a fixed number of iterations can be interpreted as a Runge-Kutta method itself. However, throughout this article, we use the acronym RKM exclusively to refer to the diagonally implicit or explicit Runge-Kutta methods against which we compare adaptive SDC.}} to reach a given order, it is more flexible due to the low order iterative solves, and attaining high order of accuracy is simple.
\deleted{For a fixed step size, SDC is also more accurate than RKM, allowing to reach the same error with fewer, larger steps.}
This can make SDC computationally competitive when medium to high accuracy is desired~\cite{FastWaveSlowWave}.
SDC has been successfully applied to complex problems that benefit from operator splitting~\cite{BOURLIOUX2003651, LAYTON2004697, minion_semi-implicit_2003} and problems where only low order solvers are available~\cite{doi:10.1137/18M1235405, p_adaptivity_phase_field}.
\added{The recent review paper by Ong and Spiteri~\cite{deferred_corrections_review} provides an overview of deferred correction methods.}

While adaptive time step selection for SDC has already been discussed in the original SDC paper~\cite{DuttSDCOriginal}, it has not yet been widely explored despite its potential to  improve computational efficiency.
Selecting time steps in SDC based on conserved quantities has been shown to mitigate order reduction~\cite{QUAIFE2016478, CiCP-24-123}.
An algorithm similar to what we propose here, based on comparing the high order SDC solution to a low order secondary solution, was shown to work well~\cite{doi:10.1137/18M1235405}.
However, both these approaches were tailored to specific problems.

This paper adopts techniques for adaptive step size selection from explicit and diagonally implicit RKM to SDC methods.
\deleted{We show that well-known results carry over to SDC, in particular improved computational efficiency and resilience against soft faults.}
What is new in this paper is that we use generic error estimates and combine a generic step size update equation known from embedded RKM with SDC.
We show for four different nonlinear problems that adaptive step size selection for SDC improves computational efficiency with little problem-specific tweaking \added{in a way that is familiar from embedded RKM}.
\added{The step size selection requires only one additional parameter compared to standard SDC, a tolerance that is used to control an estimate of the local error.}

SDC has gained popularity in the parallel-in-time (PinT) integration community due to its iterative and highly flexible nature~\cite{50yrsPinT}.
SDC-based PinT algorithms have been devised in both parallel-across-the-method~\cite{parallelSDC} as well as parallel-across-the-steps~\cite{PFASST} fashion.
We investigate adaptive step size selection for the parallel-across-the-steps Block Gau{\ss}-Seidel SDC algorithm and the parallel-across-the-method diagonal SDC approach.
So far, there are only very few studies that explore adaptive step size selection for the Parareal PinT method~\cite{MadayEtAl2020,LegollEtAl2022,KazakovEtAl2022} and none for SDC-flavored PinT algorithms.

\deleted{Although adaptivity is typically considered a means to improve computational efficiency, it can also improve resilience of simulations against faults causing silent data corruption.}
\deleted{With high-performance computing systems growing ever larger, this is becoming a concern for software developers.}
\deleted{We show that, as for RKM, error estimators in SDC used in adaptive time step selection can also be used to detect changes in data due to soft faults.}

\section{Methods and background\label{sec:methods}}
This section reviews the SDC algorithm, including its parallelizable variants diagonal SDC and Block Gau{\ss}-Seidel SDC.
It outlines adaptive step size selection strategies known from embedded RKM and proposes two ways of applying them to SDC.
\deleted{Finally, an SDC variant of the Hot Rod resilience strategy is presented.}

\subsection{Spectral Deferred Correction}
SDC methods perform numerical integration of initial value problems (IVP)
\begin{equation}
    \label{eq:ivp}
    u_t = f(u), ~~~ u(t_0) = u_0,
\end{equation}
consisting of an ODE and initial conditions, where $u$ is the solution, $f$ the right-hand side function, and the subscript $t$ indicates a derivative with respect to time.
For ease of notation we restrict the discussion to autonomous scalar problems, i.e.~$u(t)$, $u_0$, $f(u)$ $\in \mathbb{R}$ and solve a single time step from $t_0=0$ to $t_1=\Delta t$.
All derivations and results shown here can be transferred to higher dimensions, albeit with a more cumbersome notation.

We recast~\autoref{eq:ivp} in Picard form
\begin{equation}
    \label{eq:ivpPicard}
    u(t) = u_0 + \int_{0}^t f\left(u(s)\right) ds, ~~~ t \in [0, \Delta t],
\end{equation}
by integrating with respect to time.
Next, the integral is approximated with a spectral quadrature rule using $M$ collocation nodes $0\leq \tau_m\leq 1$, $m=1, ..., M$, rescaled by $\Delta t$ to cover the interval $[0, \Delta t]$.
The aim is to obtain values of the solution at the quadrature nodes and to use polynomial interpolation to approximate the continuous solution.
Using Lagrange polynomials~\cite{Interpolation}
\begin{equation}
    \label{eq:LagrangePolynomials}
    l_j^{\tau}(t)=\frac{\prod_{i=1, i\neq j}^M (t-\tau_i)}{\prod_{i=1, i\neq j}^M (\tau_j-\tau_i)}.
\end{equation}
to approximate the right-hand side function yields
\begin{equation}
    f(u(t)) \approx \sum_{j=1}^M f(u(\Delta t\tau_j)) l^{\tau}_j(t/\Delta t).
\end{equation}
By defining quadrature weights
\begin{equation}
q_{mj} = \int_{0}^{\tau_m} l_j^\tau(s) ds,
\end{equation}
we can approximate the solution at the quadrature nodes by
\begin{equation}
    \label{eq:collocation_stage}
u_m := u_0 + \Delta t\sum_{j=1}^M q_{mj}f(u(\Delta t\tau_j)) \approx u(\Delta t \tau_m).
\end{equation}
We call~\autoref{eq:collocation_stage} for $m=1, \ldots, M$ the collocation problem.
To streamline notation, we write the collocation problem in vector form
\begin{equation}
\vec{u} = \vec{u}_0 + \Delta t QF(\vec{u}),
\end{equation}
where $Q\in\mathbb{R}^{M\times M}$ is the quadrature matrix containing the weights $q_{mj}$, $F(\vec{u})=(f(u_1), ..., f(u_M))^T\in \mathbb{R}^M$ is a vector function for the temporal evolution, $\vec{u} = (u_1, ..., u_M)^T \in \mathbb{R}^M$ is the vector carrying the approximate solutions at the quadrature nodes and $\vec{u}_0=(u_0, ..., u_0)^T \in \mathbb{R}^M$ the initial conditions.

This collocation problem corresponds to a fully implicit Runge-Kutta method.
It can be solved directly, but doing so is expensive because $Q$ is dense and the equations for the stages~(\autoref{eq:collocation_stage}) are all coupled.
Certain types of nodes, such as Gau{\ss}-Legendre, Gau{\ss}-Lobatto and Gau{\ss}-Radau can achieve super-convergence up to order $2M$ at the right boundary with only $M$ nodes~\cite[Theorems 5.2, 5.3 and 5.5]{hairer_wanner_II}.

Using Picard iteration to solve the collocation problem provides a simple iterative scheme
\begin{align}
\vec{u}^{k+1} &= \vec{u}^{k} + \vec{r}^k\\
\vec{r}^{k} &:= \vec{u}_0 + \Delta tQF(\vec{u}^k) - \vec{u}^k
\end{align}
where the solution in iteration $k$ is improved by adding the residual $\vec{r}^k$.
However, this method only converges for small $\Delta t$ or (very) non-stiff ODEs.
SDC preconditions Picard iterations with a low order method to achieve convergence also for large time step size.
\added{As we will now illustrate, this can be seen as integrating an equation for the error with said low order method.}

The quadrature rule integrates the polynomial approximation exactly, i.e.~$\int_{0}^{t} F(\vec{u})\vec{l}^{\tau}(s/\Delta t)ds=\Delta t(QF(\vec{u}))\vec{l}^{\tau}(t/\Delta t)$.
Thus, if the error of the current polynomial approximation $\delta^k(t)=u(t)-\vec{u}^k_m\vec{l}^{\tau}(t/\Delta t)$ is plugged into~\autoref{eq:ivpPicard}, we get
\begin{equation}
    \delta^{k}(t) - \int_{0}^t\left(f\left(\vec{u}^k\vec{l}^{\tau}(s/\Delta t)+\delta^{k}(s)\right) - F\left(\vec{u}^k\right)\vec{l}^{\tau}(s/\Delta t)\right)ds = \vec{r}^k\vec{l}^{\tau}(t/\Delta t).
\end{equation}
The integral is approximated by a simpler quadrature rule $Q_\Delta$, typically referred to as preconditioner, and the resulting nonlinear system 
\begin{equation}
    \vec{\delta}^k - \Delta t Q_\Delta\left(F(\vec{u}^k + \vec{\delta}^k) - F(\vec{u}^k)\right) = \vec{r}^k.
\end{equation}
has to be solved in each iteration.
The solution is then updated by adding the correction, $\vec{u}^{k+1} = \vec{u}^k + \vec{\delta}^k$.
By expanding the residual and plugging in the refinement equation, we can eliminate the defect and simplify the SDC iteration to
\begin{equation}
    \label{eq:sdcIteration}
    \left(I_M - \Delta t Q_\Delta F\right)\left(\vec{u}^{k+1}\right) = \vec{u}_0 + \Delta t\left(Q - Q_\Delta\right)F\left(\vec{u}^k\right),
\end{equation}
with $I_M$ the $M\times M$ identity matrix.

The preconditioner $Q_\Delta$ is typically chosen to be a lower triangular matrix such that the system can be solved with forward substitution.
In the context of partial differential equations with $N$ degrees of freedom, the collocation problem is a $NM\times NM$ system and iterating $K$ times with forward substitution allows to instead solve $KM$-many $N\times N$ systems.
Since the algorithm proceeds from one line of the system to the next, SDC iterations are often referred to as ``sweeps''.

\paragraph{Preconditioners and diagonal SDC}
In the original derivation of SDC, implicit or explicit Euler were proposed for solving the error equations~\cite{DuttSDCOriginal}.
These are first order quadrature rules, integrating from node to node.
This means they increase the order of accuracy by one up to the order of the underlying collocation problem, provided they converge at all.
The preconditioner corresponding to implicit Euler, for example, reads:
\begin{equation}
    Q_\Delta^\mathrm{IE} = \left(\begin{matrix}
		\tau_2 - \tau_1 & 0 & 0 & \hdots & 0 \\
		\tau_2 - \tau_1 & \tau_3 - \tau_2 & 0 & \hdots & 0 \\
		\vdots               & \vdots                & \ddots & \hdots & 0 \\
		\tau_2 - \tau_1 & \tau_3 - \tau_2 & \hdots &  \hdots & \tau_{M} - \tau_{M-1} \\
		
	\end{matrix}\right).
\end{equation}
Higher order preconditioners such as RKM can sometimes increase the order of accuracy by more than one with each iteration~\cite{highOrderPrecondsEquidistant}.
However, lack of smoothness in the error can limit gains to one order per iteration regardless of the order of the preconditioner, particularly when non-equidistant nodes are used~\cite{high_order_preconditioner}.

Other interpretations of SDC do not rely on the preconditioner being consistent with an integration rule~\cite{SDC_convergence}.
The LU preconditioner~\cite{LU}, where $Q_{\Delta} = U^\mathrm{T}$, with $LU = Q^\mathrm{T}$, is very effective for stiff problems.
Be aware that most of these preconditioners cannot be interpreted as A-stable time marching schemes and usually pose some restrictions on the domain of convergence of the collocation problem.
\deleted{In the case of the LU preconditioner, a gain of exactly one order per sweep is not always observed for any given step size, see also the discussion in kremling2021convergence.
For this reason, we stick to the implicit Euler preconditioner for most simulations.}

Another promising class of preconditioners use a diagonal $Q_{\Delta}$ and allow updating all nodes in parallel~\cite{parallelSDC}.
This is a small-scale PinT algorithm since the number of nodes is typically small, but it can be combined with other PinT algorithms that solve multiple steps concurrently.
Numerical experiments suggest that best parallel efficiency is obtained when maximizing parallelism across-the-nodes and spending the remaining computational resources on parallelizing across-the-steps~\cite{schobel_pfasst-er_2020}.

Good diagonal preconditioners can be derived, for instance, by minimizing the spectral radius of the SDC iteration matrix in the stiff or non-stiff limit of the test equation.
When running diagonal SDC, we use preconditioners derived in this way\deleted{~from}~\cite{GT_precons}, namely \deleted{MIN-SR-NS for problems with moderate stiffness and }MIN-SR-S\added{, which is suitable} for stiff problems.
The performance of diagonal SDC variants combined with space-time parallelization \replace{is}{has also been} investigated\deleted{in}~\cite{parallel_sdc_hpc}.

\added{Since the LU and diagonal preconditioner cannot guarantee an increase by exactly one order per sweep~\cite{kremling2021convergence,GT_precons}, we stick to the implicit Euler preconditioner for most simulations.}

\paragraph{Inexact SDC iterations}
Performing the implicit solves inside the SDC iteration to full accuracy often comes at no benefit to overall accuracy because the preconditioner itself is only approximate.
Reducing tolerances or strict limits on the number of allowed iterations for the implicit solver can be used to avoid over-solving and to improve overall computational efficiency of SDC~\cite{inexactSDC}.
Optimal tolerances can be derived, but require realistic work and error models~\cite{weiser_inexactness_2018}.
Still, efficiency can be gained even with sub-optimal tolerances when adjusting the tolerance of an inner solver based on the outer residual or when allowing only few iterations of an inner solver.
Because SDC provides good initial guesses for the nonlinear solver, this can lead to very efficient schemes~\cite{inexactSDC}.

For simplicity, we fix the ratio of inner tolerance to outer residual based on heuristics to \replace{$1/10$}{some fixed value}.
Additionally, we can only employ inexactness if knowing the order after every iteration is not required because otherwise we cannot guarantee that the inner solver is accurate enough to increase the order by one.

\paragraph{Implicit-explicit splitting}
When solving problems with stiff and non-stiff components it is possible to treat only the stiff part implicitly and the remainder explicitly.
This is called implicit-explicit (IMEX) splitting and can easily be used in SDC~\cite{minion_semi-implicit_2003, FastWaveSlowWave}.
The IMEX-SDC iteration reads
\begin{align*}
    (1-\Delta t \Tilde{q}_{m+1,m+1}^I f^I)(u_{m+1}^{k+1}) = u_0 & + \Delta t \sum_{j=1}^m \left(\Tilde{q}_{m+1, j}^I f^I + \Tilde{q}_{m+1,j}^E f^E\right)(u_{j}^{k+1}) \\ 
    & -\Delta t \sum_{j=1}^{m+1} \left(\Tilde{q}_{m+1,j}^I f^I + \Tilde{q}_{m+1,j}^E f^E\right)(u_{j}^{k}) \\
    & + \Delta t \sum_{j=1}^M q_{m+1, j} (f^I + f^E)(u_{j}^k).
\end{align*}
with superscripts $I$ and $E$ referring to the implicit and explicit part and $\Tilde{q}$ being the entries of the preconditioners.
The scheme typically has the same order of accuracy as the non-split version, although order reduction may occur~\cite{minion_semi-implicit_2003, IMEXFACEOFF}.
IMEX-SDC has been shown to outperform DIRK-based IMEX Runge-Kutta methods for incompressible flow simulations in wall-time measurements~\cite{IMEXFACEOFF}.

\paragraph{Dense output}
As the solution of the collocation problem is an approximation by a polynomial of degree $M$, a natural continuous extension is suggested by evaluating this polynomial anywhere within the interval.
We refer to this as ``dense output'' property of the collocation problem~\cite[Sect. II.6]{Hairer_Wanner}.
Keep in mind that, while the solution at the boundary is up to order $2M$ accurate for $M$ nodes, the accuracy inside the interval is order $M$.

\subsection{Adaptive step size selection for SDC}
\label{subsec:adaptive_time-stepping}
Adaptive selection of the time step size is useful when the rate of change of the solution is not uniform across the computational domain.
We transfer well-known concepts from embedded Runge-Kutta methods~\cite{Hairer_Wanner}  to SDC.
First, we need to estimate the local error, which we then aim to control by choosing an appropriate step size.

The error estimation works by computing two solutions to the same initial value problem with a different order of accuracy.
The difference between the solutions is a reasonable estimate of the local error of the less accurate method
\begin{align}
\label{eq:error_estimate}
    \begin{split}
    \epsilon &= \|u^{(p)} - u^{(q)}\|_{\new{\infty}} \\
    &= \|(u^{(p)} - u^{*})-(u^{(q)} - u^{*})\|_{\new{\infty}} \\
    &= \|\delta^{(p)} - \delta^{(q)}\|_{\new{\infty}} = \new{\|}\delta^{(p)}\new{\|_\infty} + \mathcal{O}(\Delta t^{q+1}),
    \end{split}
\end{align}
where $u^{(p)}$, $u^{(q)}$ are the solutions \new{at the end of the time interval,} obtained by integration schemes of order $p$ and $q$ with $q > p$. $u^*$ marks the exact solution, $\delta$ denotes the local error with analogous meaning of the superscript and $\epsilon$ is the estimate of the local error.
Once $\epsilon$ is known, an optimal step size
\begin{equation}
     \Delta t_\mathrm{opt} = \beta \Delta t \left(\frac{\epsilon_\mathrm{TOL}}{\epsilon}\right)^{1/(p+1)}
     \label{eq:step_size_update}
\end{equation}
can be estimated such that $\epsilon\approx\epsilon_\mathrm{TOL}$.
Here, $\epsilon_\mathrm{TOL}$ is the user-defined tolerance for the local error and $\beta$ is a safety factor, usually $\beta=0.9$.
This update equation is based on the order of accuracy $p$ of the time-marching scheme~\cite{Hairer_Wanner, RK54}.
As is also common in embedded RKM, we use ``local extrapolation'', meaning we advance using the higher order solution, even though we control the error of the lower order one.
Crucially, we check if the local error estimate falls below the desired accuracy and we move on to the next step with $\Delta t_\mathrm{opt}$ only if it does.
If we fail to satisfy the accuracy requirements, we recompute the current step with $\Delta t_\mathrm{opt}$.
While the time-scale of the problem may have changed in the next step, heuristically, $\Delta t_\mathrm{opt}$ often appears to be a good guess for the optimal step size for the next step as well.

\paragraph{Adaptive selection of $\Delta t$}
Since the order increases by one with each SDC iteration (up \new{to} the order of the underlying collocation problem), the increment can be used directly as the error estimate.
While adaptivity based on the increment was already proposed in the original SDC paper~\cite{DuttSDCOriginal}, they employ a simpler step size update equation based only on doubling or halving.
\added{A similar approach, although not yet called SDC, was studied by van der Houwen and Sommeijer~\cite{vanderHouwen1992embedded, VANDERHOUWEN1990111} for other preconditioners.}

Algorithm~\ref{alg:dtadaptivity} shows in pseudo-code the algorithm resulting from combining an increment based error estimation with~\autoref{eq:step_size_update}.
\begin{algorithm}[t]
\caption{SDC with $\Delta t$-adaptivity}\label{alg:dt_adaptivity}
\begin{algorithmic}
\State $u^0 \gets u_0$ 
\State $k \gets 1$
\While{$k\leq k_\mathrm{max}$}
   \State $u^k \gets \mathrm{SDC\ iteration\ applied\ to\ } u^{k-1}$
   \State $k \gets k+1$
\EndWhile

\State $\epsilon \gets \|u^{k_\mathrm{max}} - u^{k_\mathrm{max}-1}\|$
\State $\Delta t \gets  \beta \Delta t \left(\frac{\epsilon_\mathrm{TOL}}{\epsilon}\right)^{1/k_\mathrm{max}}$

\If{$\epsilon\gt\epsilon_\mathrm{TOL}$}
    \State Restart current step with $u_0$
\Else
    \State Move on to next step with $u^{k_\mathrm{max}}$
\EndIf
\end{algorithmic}
\label{alg:dtadaptivity}
\end{algorithm}
Since this approach modifies only $\Delta t$ but keeps the number of iterations constant, we refer to it as $\Delta t$-adaptivity.
Order reduction may be observed for very stiff problems, requiring some extra care~\cite{HUANG2006633, layton_implications_2005}.

\paragraph{Adaptive selection of $\Delta t$ and $k$}
By using the dense output property of the converged collocation problem, we can design an approach to choose both the time step $\Delta t$ and the number of iterations $k$ adaptively.
For the fully converged collocation solution, the $M+1$ values $u_0, u_1, \ldots, u_M$ \added{at the collocation points $\tau = \left\{ \tau_i : i = 0, \ldots, M \right\}$} define a polynomial on $[0, \Delta t]$ that provides an  $M$-th order accurate approximation at any point $t \in [0, \Delta t]$.
By removing \replace{one}{the} collocation point \added{$\tau_{M-1}$ and using a set of nodes
\begin{equation}
    \tau^* = \{\tau_i : i=0, \ldots, M-2, M\},
\end{equation}
}we can construct \added{a polynomial of degree $M-1$ and evaluate it at $\tau_{M-1}$ to produce} an $M-1$-st order accurate approximation \replace{anywhere within the interval}{
\begin{equation}
    u^{(M-1)}_{M-1} = \sum_{i=0}^{M-2} u^{k}_i l^{\tau^*}_i(\tau_{M-1}) + u^{k}_M l^{\tau^*}_{M-1}(\tau_{M-1})    
\end{equation}
of the solution at the removed point.
Here, $l^{\tau^*}_i$ are the Lagrange polynomials of the $i$-th node in the set $\tau^*$.
We can then use the difference 
\begin{align}
    \epsilon &= \|u^{(M-1)}_{M-1} - u^k_{M-1}\|_\infty,
\end{align}
as an error estimate.}
\deleted{Taking the difference this generates an error estimate of order $M$.}
\new{Tests not documented here suggest that the node at which to compute the error can be chosen more or less arbitrarily. 
We choose node $M-1$ in our examples since this performed well and leave a more detailed mathematical investigation for future work.}

For very large step sizes the SDC iteration may not converge with arbitrary preconditioner.
To mitigate, we introduce a relative limit $\gamma=4$ on how much the step size is allowed to increase and set $k_\mathrm{max} = 16$\reviewcomment{was 20 before by mistake. In the experiments we had always used 16.}, so that we can return to the last step size that allowed convergence if the desired residual tolerance was not achieved after $k_\mathrm{max}$ iterations.

Since this approach does not require a fixed $k$, we can simply keep iterating until convergence and thus let SDC choose the number of iterations as well as the step size.
We therefore refer to this approach as $\Delta t$-$k$-adaptivity.
Algorithm~\ref{alg:dtkadaptivity} sketches the algorithm in pseudo-code.
\begin{algorithm}
\caption{SDC with $\Delta t$-$k$-adaptivity}
\begin{algorithmic}
\State $\vec{u}^0 \gets \vec{u}_0$ 
\State $k \gets 1$
\State $r \gets \|\vec{u}_0 +\Delta t QF(\vec{u}^0)- \vec{u}^0\|_{\new{\infty}}$

\State $r_\mathrm{prev}\gets \infty$
\State $r_\mathrm{max}\gets 10^9$
\State $no\_convergence\gets False$
\While{$r\gt r_\mathrm{tol}$ and not $no\_convergence$}
   \State $\vec{u}^k \gets \mathrm{inexact\ SDC\ iteration\ applied\ to\ } \vec{u}^{k-1}$
   \State $r \gets \|\vec{u}_0 +\Delta t QF(\vec{u}^k)- \vec{u}^k\|_{\new{\infty}}$
   \If{$r > r_\mathrm{max} \mathrm{~or~} r > r_\mathrm{prev} \mathrm{~or~} k = k_\mathrm{max}$}
        \State $no\_convergence \gets True$
   \EndIf
   
   \State $r_\mathrm{prev}\gets r$
   \State $k \gets k+1$
\EndWhile

\If{$no\_convergence$}
    \State $\Delta t \gets \Delta t / \gamma$
    \State Restart current step with $u_0$
\Else

\State $\tau^* \gets \{\tau_i, i\neq M-1\}$
\State $\epsilon \gets \|\sum_{i=0}^{M-2} u^{k}_i l^{\tau^*}_i(\tau_{M-1}) + u^{k}_M l^{\tau^*}_{M-1}(\tau_{M-1}) - u^{k}_{M-1}\|_{\new{\infty}}$

\State $\Delta t \gets  \max \left(\gamma,  \beta \left(\frac{\epsilon_\mathrm{TOL}}{\epsilon}\right)^{1/M}\right) \Delta t $

\If{$\epsilon\gt\epsilon_\mathrm{TOL}$}
    \State Restart current step with $u_0$
\Else
    \State Move on to next step with $u^{k}$
\EndIf
\EndIf
\end{algorithmic}
\label{alg:dtkadaptivity}
\end{algorithm}
\reviewcomment{Removed brackets around $\tau$ in Algorithm \autoref{alg:dtkadaptivity}.}
\new{Note that the residual tolerance controls how accurately we solve the collocation problem while the error tolerance controls how close we are to the continuous solution of the ODE.
Therefore, some care must be taken to reasonably balance the two - a very small residual with a large error tolerance, for example, will lead to a very good approximation of the collocation polynomial that is a poor approximation of the continuous solution.
In our numerical examples, we choose a residual tolerance $r_\mathrm{TOL}$ a few orders of magnitude smaller than $\epsilon_\mathrm{TOL}$.}

\paragraph{Mitigating the cost of restarts}
Restarting steps from scratch in Algorithms~\autoref{alg:dt_adaptivity} and~\autoref{alg:dtkadaptivity} is expensive but SDC offers a unique way to reduce this overhead:
Using the dense output property, we can evaluate the polynomial of a restarted step at the new collocation nodes resulting from shorter step size, thus re-using the previously computed solution.
Despite being not accurate enough to satisfy the prescribed tolerance, this is a good initial guess, and the SDC iteration is expected to converge in fewer iterations.
This is sensible to do in the $\Delta t$-$k$-adaptivity strategy, but only when the collocation problem of the restarted step has converged, i.e., the SDC residual is small enough, as otherwise the solution is not a useful approximation. 
However, with the $\Delta t$-adaptivity strategy, a gain in efficiency is unlikely as the step size update equation will overestimate the optimal step size due to a one-time exceptionally good initial guess which the next step will not have access to.

\paragraph{Advantages of $\Delta t$-and $\Delta t$-$k$-adaptivity over $k$-adaptivity}
Even without local error estimation, SDC can be used adaptively with a fixed step size by simply iterating until the SDC residual is below a set tolerance.
We call this approach $k$-adaptivity.
Our newly introduced $\Delta t$- and $\Delta t$-$k$-adaptivity approaches have two major advantages over $k$-adaptivity.
First, the step size is a floating point number, which can be finely adjusted, whereas the number of iterations is an integer allowing the scheme\deleted{s} much less control.
Second, the local error estimate, \autoref{eq:error_estimate}, measures the error with respect to the continuous solution.
By contrast, the SDC residual is only a measure of the iteration error with respect to the discrete collocation solution.
While a low residual indicates that the collocation problem has been solved to high accuracy, the truncation error of the method might still be large.

\subsection{Pipeline-based parallelism: Block Gau{\ss}-Seidel SDC}
\label{sec:GSSDC}
A relatively straightforward pipeline-based parallel-in-time variant of SDC can be constructed by solving multiple time steps simultaneously in block Gau{\ss}-Seidel fashion~\cite{GSSDCfirst}:
Instead of waiting for the previous step to converge to full accuracy before starting a new step, we begin solving a new step as soon as a single iteration has been performed on the previous step in the block and keep refining the initial conditions with the iterates from the previous step between iterations.
\begin{figure}
    \centering
    \includegraphics[width=0.5\linewidth]{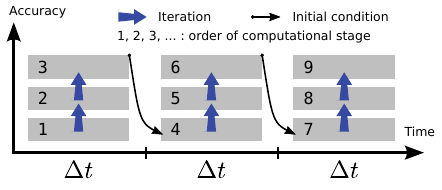}%
    \includegraphics[width=0.5\linewidth]{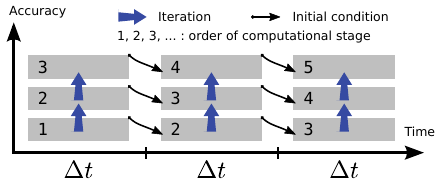}
    \caption{Left: Sequential SDC. Right: Block Gau{\ss}-Seidel SDC. Sequential time stepping uses the converged solution of the last step as initial conditions, whereas the parallel version receives refined initial conditions between iterations.}
    \label{fig:bgs-sdc}
\end{figure}
\autoref{fig:bgs-sdc}~illustrates both sequential SDC one the left and parallel-across-the-steps Block Gau{\ss}-Seidel SDC (GSSDC) on the right.
A very similar parallel-in-time (PinT) algorithm based on pipelining more general deferred corrections, Revisionist Integral Deferred Correction\deleted{s} (RIDC) \cite{RIDC}, has been shown to provide good speedup.
We focus on GSSDC as it is at the heart of the large-scale PinT algorithm PFASST~\cite{PFASST}.
In contrast to Guibert et al.~\cite{GSSDCfirst}, however, we do not use any overlap between the collocation nodes of different time steps.
\new{When using adaptive step size selection, we use a single $\Delta t$ for all steps in the block.}

It has been demonstrated that the convergence order is maintained when increasing the number of steps $N$ in a block of GSSDC~\cite{thibaut_unified_2022}.
In particular, it is the same as for single step SDC, which means we can employ the adaptive step size controller with no modification.
However, the results will not be identical to serial SDC due to the inexactness in the initial conditions.

While we can still estimate the local error by means of the increment, its interpretation is different in the multi-step version.
Within a block, we are solving all steps with a first order accurate method, then we solve all steps with a second order accurate method and so on.
This means that the increment is an estimate of the global error within the block.
As the whole problem is divided into multiple blocks, we should still view this error as local in the context of the global time domain, but we need to be aware that the same local tolerance $\epsilon_\mathrm{tol}$ applied to a block of multiple steps should result in smaller local errors in each step inside the block.

When the error estimate exceeds the prescribed tolerance, we have two choices for restarting: We can restart from the initial conditions of the first step in the block, or we can restart from the initial conditions of the first step in the block where the error estimate exceeds the threshold.
The first is the more rigorous strategy while the second strategy is heuristic but can improve performance at a slight cost to accuracy.
In the tests included here, we show only the second strategy because it performed consistently better.

\reviewcomment{deleted introduction for resilience section}

\section{Benchmark problems}\label{sec:Problems}
We use one nonlinear ODE and three nonlinear PDEs as benchmarks to test the performance of our adaptive strategies.
Each benchmark problem poses different challenges for numerical time-integration.
All implementations are publicly available on~\href{GitHub}{https://github.com/Parallel-in-Time/pySDC}\footnote{\url{https://github.com/Parallel-in-Time/pySDC}} as part of the pySDC library~\cite{pySDC}.
Wall-clock times are measured on the JUSUF supercomputer~\cite{JUSUF} at Forschungszentrum J\"ulich.

\subsection{Van der Pol\label{sec:vdPProblem}}
The van der Pol equation
 \begin{align}
     u_{tt} - \mu \left(1-u^2\right)u_t + u = 0, \\
     u(t=0) = u_0, \hspace{0.2cm}
     u_t(t=0) = u_0',
 \end{align}
is named after a Dutch electrical engineer who used the equation to study the behavior of vacuum tubes in radios~\cite{van_der_pol_1926}.
Here, $\mu$ is a parameter controlling the nonlinearity, $u$ is the solution and the subscript $t$ marks a derivative with respect to time.
In our tests, we set \replace{$u_0 = 2$}{$u_0 = 1.1$} and $u_0' = 0$, \replace{$\mu=5$}{$\mu=1000$} and solve up to \replace{$t=11.5$}{$t=20$}.
In the pySDC implementation, we introduce $v(t) = u_t(t)$, rewrite the van der Pol equation as a first order system and use a Newton scheme to solve the nonlinear systems within the SDC sweeps.

\begin{figure}
    \centering
    \includegraphics{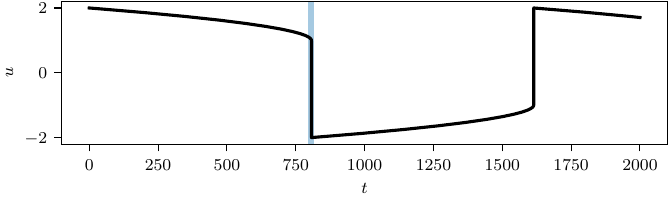}
    \caption{Solution of a van der Pol problem for \replace{$\mu=5$}{$\mu=1000$} over time. The solution is oscillating on two time-scales. In order not to over-resolve the slow parts, the resolution has to be adjusted during runtime. \added{In numerical tests, we solve only the shaded transition at $t\approx 800$ as this is already very expensive with fixed step size schemes. Note that at this high value of $\mu$, the problem is extremely stiff.}}
    \label{fig:vdpSolution}
\end{figure}

For $\mu=0$, we recover the harmonic oscillator, but with increasing $\mu$ the problem becomes increasingly stiff.
Van der Pol describes the problem with $\|\mu\| \ll 1$ as modelling free oscillations of a triode oscillator, whereas $\|\mu\| \gg 1$ models a free relaxation oscillation~\cite{vdp_modelling}.
This is a useful test problem for adaptive step size control as the nonlinear damping introduces a second time-scale to the oscillation, see~\autoref{fig:vdpSolution}.
\deleted{The pattern of oscillations in the van der Pol equation is mostly determined by the $\mu$ parameter and not by the initial condition.
This means, in particular, that the problem is not overly sensitive to perturbations, making it rather forgiving towards soft faults.}

We use the SciPy~\cite{SciPy} method \textsc{solve\_ivp} from the \textsc{integrate} package with an explicit embedded Runge-Kutta method of order 5(4)~\cite{DORMAND198019} with tolerances close to machine precision to obtain reference solutions.
\added{When using $\Delta t$-$k$-adaptivity, we select $r_\mathrm{TOL}=10^{-5}\epsilon_\mathrm{TOL}$, and stop the Newton scheme at a tolerance of $10^{-5}r$, with $r$ the current SDC residual, or after a maximum of 9 iterations.}

\subsection{Quench \label{sec:QuenchProblem}}
This is a simplified model of temperature leaks in superconducting magnets provided by Erik Schnaubelt~\cite[Section 4.3]{Schnaubelt_quench}.
Once the temperature exceeds a certain threshold, superconductivity ceases and runaway heating of the magnet sets in.
This effect has led to the explosion of large magnets at the Large Hadron Collider~\cite{LHCexplosion}.

\begin{figure}
    \centering
    \includegraphics{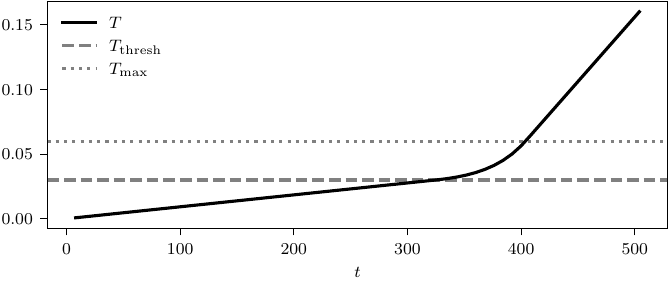}
    \caption{Solution of the Quench problem. Shown is the maximal temperature $T(t) = \|\new{u}(x,t)\|_{L^\infty(\Omega)}$ across the spatial domain over time. We see a slow heating up to about $t=320$, after which the temperature in some parts of the domain exceeds the threshold value and a linear transition towards runaway heating is entered. Physically, this means the magnet stops being superconducting, which can have catastrophic effects in particle accelerators.}
    \label{fig:quenchSolution}
\end{figure}

The model consists of a one-dimensional heat equation with a nonlinear source term heating parts of the domain.
For the boundary conditions, we choose Neumann-zero to treat the magnet as completely isolated from the environment, except for the leak.
Due to superconductivity, the diffusivity is high, making the problem very stiff and prohibiting the use of explicit time-stepping schemes.
The equation reads
\begin{align}
    C_V u_t  - \kappa \Delta u = & Q(u), &\\
    Q(u) = & Q_\mathrm{max}\times \begin{cases}
        1, & x \in (0.45, 0.55), \\
        f(u), & \mathrm{else},
    \end{cases} \\
    f(u) = & \begin{cases}
        0, & u < T_\mathrm{thresh}, \\
        \frac{u - T_\mathrm{thresh}}{T_\mathrm{max} - T_\mathrm{thresh}}, & T_\mathrm{thresh} \leq u < T_\mathrm{max}, \\
        1, & T_\mathrm{max} \leq u,
    \end{cases} \\
    \Omega \in & ]0, 1[, \\
    u_x = & 0,~x \in \partial \Omega, \\
    u(t=0) = & 0,
\end{align}
with the Laplacian $\Delta$ and parameters
\begin{align*}
    C_V = & 1000, \\
    \kappa = & 1000, \\
    T_\mathrm{thresh} = & 10^{-2}, \\
    T_\mathrm{max} = & 2\times 10^{-2}, \\
    Q_\mathrm{max} = & 1.
\end{align*}
We solve until $t=500$ and use a Newton scheme for implicit solves.
\autoref{fig:quenchSolution} shows the solution over time.
The behavior is fairly simple, except for the transition from superconductivity to runaway heating, which is challenging to resolve accurately.

We again rely on the SciPy method \textsc{solve\_ivp} to generate reference solutions, but, as the problem is very stiff, use an implicit backward differentiation formula~\cite{BDF}.
\added{When using $\Delta t$-$k$-adaptivity, we select $r_\mathrm{TOL}=10^{-1}\epsilon_\mathrm{TOL}$, and stop the Newton scheme at a tolerance of $10^{-1}r$ or after a maximum of 5 iterations.}

\subsection{Nonlinear Schr{\"o}dinger \label{sec:SchroedingerProblem}}
\begin{figure}
    \centering
    \includegraphics{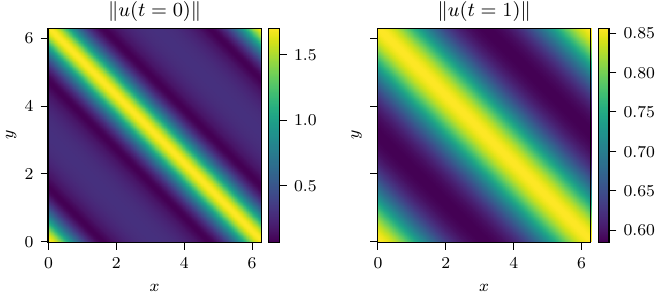}
    \caption{\reviewcomment{This figure is new}\added{Plot of the initial condition for the Schr{\"o}dinger example (left) and the solution at the end of the interval under consideration (right). The initial conditions are purely real, but the solution shifts to the complex domain as the simulation progresses.}}
    \label{fig:SchroedingerSolution}
\end{figure}

The focusing nonlinear Schr{\"o}dinger equation is a wave-type equation that describes problems such as signal propagation in optical fibers~\cite{NLS_stuff}.
The formulation we solve can be written as
\begin{align}
    \label{eq:schroedinger}
    u_t =& i\Delta u + 2i \|u\|^2u,\\
    u(t=0) =& \frac{1}{\sqrt{2}} \left(\frac{1}{1 - \cos(x + y)/ \sqrt{2}} - 1\right),
\end{align}
\reviewcomment{fixed sign error in \autoref{eq:schroedinger}}
with $i$ being the imaginary unit.
We solve~\autoref{eq:schroedinger} on a two-dimensional spatial domain with periodicity $2\pi$ up to $t=1$ using fast Fourier transforms.
We use implicit-explicit (IMEX) splitting to integrate the Laplacian implicitly and the nonlinear term explicitly.
The global error is computed with respect to the analytic solution ~\added{\cite[Equation (39)]{Schroedinger_exact_original}}\reviewcomment{Switched reference for exact solution.}
\added{See \autoref{fig:SchroedingerSolution} for a visualization of the solution.}
\added{When using $\Delta t$-$k$-adaptivity, we select $r_\mathrm{TOL} = 10^{-4}\epsilon_\mathrm{TOL}$}

\subsection{Allen-Cahn}
\begin{figure}
    \centering
    \includegraphics{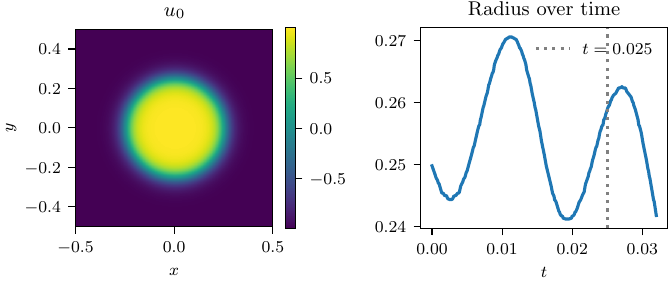}
    \caption{Left: Initial condition consisting of a circle of high phase embedded in the low phase. Right: \replace{Decrease}{Evolution} of the radius of the circle over time. The simulation is terminated at $t=0.025$\deleted{before the circle has disappeared}.}
    \label{fig:ACSolution}
\end{figure}

The Allen-Cahn variant considered here is a two-dimensional reaction-diffusion equation with periodic boundary conditions \added{that can be used to model transitions between two phases}.
\replace{We choose initial conditions that cause a phase transition in the pattern of a contracting circle~cite{AC}.}{We choose initial conditions representing a circle of one phase embedded in the other phase.
We add time-dependent forcing, such that the circle alternates between growing and shrinking.}
\begin{align}
    u_t &= \Delta u - \frac{2}{\epsilon^2}u(1-u)(1-2u) - 6u(1-u)f(u, t),\\
    f(u, t) &= \frac{\sum \left(\Delta u - \frac{2}{\epsilon^2}u(1-u)(1-2u)\right)}{\sum 6u(1-u)}\left(1 - \sin{\left(4\pi \frac{t}{0.032}\right)}\times 10^{-2}\right)\\
    u_0(x) &= \tanh{\left(\frac{R_0\|x\|}{\sqrt{2}\epsilon}\right)},\\
    x &\in {[-0.5, 0.5[^2}\\
    \epsilon &= 0.04,\\
    R_0 &= 0.25,
\end{align}
\reviewcomment{Changed problem definition.}
\deleted{In the sharp interface limit $\epsilon\rightarrow 0$, the radius of the circle shrinks as $r(t)=\sqrt{R^2 - 2t}$.}
\autoref{fig:ACSolution} shows the initial condition \added{and the evolution of the radius over time}.
\deleted{The IMEX Euler solver inside the SDC iterations is not able to resolve the front for $\Delta t>\epsilon^2$~cite{AC}.
Including a safety factor of $0.8$, we limit the step size to $\Delta t\leq0.8\epsilon^2$.}
Similar to our approach for the Schr{\"o}dinger problem, we use IMEX splitting, integrating the Laplacian discretized with a spectral method implicitly while treating the nonlinear term explicitly.
\added{The IMEX scheme is not unconditionally stable, providing an upper limit on the step size and adding to the challenges for adaptive step size selection.}
We again compute the error with respect to the SciPy method \textsc{solve\_ivp}.
\added{When using $\Delta t$-$k$-adaptivity, we select $r_\mathrm{TOL} = 10^{-3}\epsilon_\mathrm{TOL}$}

\deleted{Since there are no major changes in time-scale in this problem, adaptivity will have minor benefits at best.
Instead, we use this benchmark to confirm that adaptivity does not add substantial overhead.}
\section{Numerical results}
\label{sec:Results}

We investigate performance of the methods from \autoref{sec:methods} applied to the problems from \autoref{sec:Problems} with respect to computational efficiency \replace{and resilience.
We first consider SDC alone, then compare the results to RKM and investigate resilience at the end.}{and then compare them to established RKM.}

\subsection{Computational efficiency\label{sec:Efficiency}}
\begin{figure}
    \centering
    \includegraphics{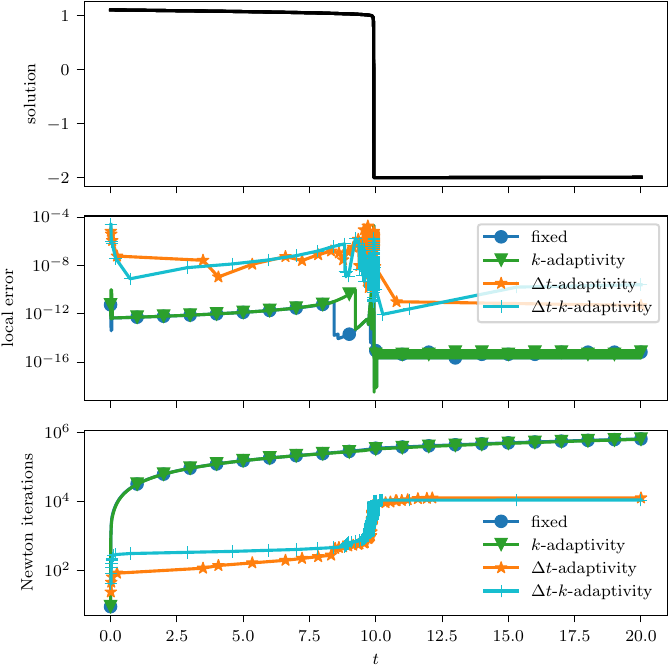}
    \caption{\reviewcomment{Swapped figure, had $\mu=5$ before, now $\mu=1000$.}SDC for the van der Pol problem \added{with $\mu=1000$, solved} with four different strategies. The top panel shows the solution, the middle panel shows the resolution via the local error compared to a reference solution and the bottom panel shows how many Newton iterations are needed to reach the respective time, which is a good indicator of computational cost. \replace{The fixed scheme with constant $k$ and $\Delta t$ requires the most iterations and the local error varies by more than four orders of magnitude. Adaptively choosing the iteration number leads to less pronounced over-resolving of slow parts, but using $\Delta t$-adaptivity to fine-tune the step size to the problem's time-scale allows to save even more iterations without sacrificing accuracy in terms of local error. The $\Delta t$-$k$-adaptive scheme shows greater variability in accuracy, but maintains the same minimum accuracy at no additional cost. While it does select smaller step sizes, fewer Newton iterations are performed per step.}{We show only a fast transition which requires very small step sizes due to extreme stiffness. The solution is periodic with the next transition appearing at around $t=600$ (see \autoref{fig:vdpSolution}). Even in this short time period, schemes with fixed step size need about 70 times as many Newton iterations as adaptive ones, with the difference only becoming larger during the long period of little action until the next transition. Note that each marker represents a single step for methods with adaptive $\Delta t$ or 10000 for the fixed or $k$-adaptive methods.}}
    \label{fig:vdpEfficency}
\end{figure}

We demonstrate the adaptive resolution capabilities in detail for the van der Pol problem first, as it easy to visualize.
\autoref{fig:vdpEfficency} shows the solution (upper), local error (middle) and computational work, measured in total number of required Newton iterations (lower) for the van der Pol equation.
\deleted{SDC with a fixed time step over-resolves significantly in many areas, delivering a local error of $10^{-11}$, much smaller than the maximum local error of slightly above $10^{-7}$.
The adaptive variants of SDC over-resolve much less, which is reflected in the significantly lower number of required Newton iterations.
Both $\Delta t$- and $\Delta t$-$k$-adaptive SDC perform very similarly, requiring only around $5000$ instead of $15000$ Newton iterations over the course of the simulation.}
\added{Resolving the transitions requires a very small step size due to the high stiffness of the problem, but in between transitions, a much larger step size suffices.
We select parameters such the maximal local error during the transition is on the order of $10^{-5}$.
For fixed step size schemes, the resulting step size is so small, that we solve to machine precision outside of the transition.
The associated computational effort to cover this short simulation time is approximately 70 times of what is required by step size adaptive strategies to meet the same accuracy requirements.
In tests not shown here, we found similar, although less pronounced, trends for smaller values of $\mu$.
}

\begin{figure}
    \centering
    \includegraphics{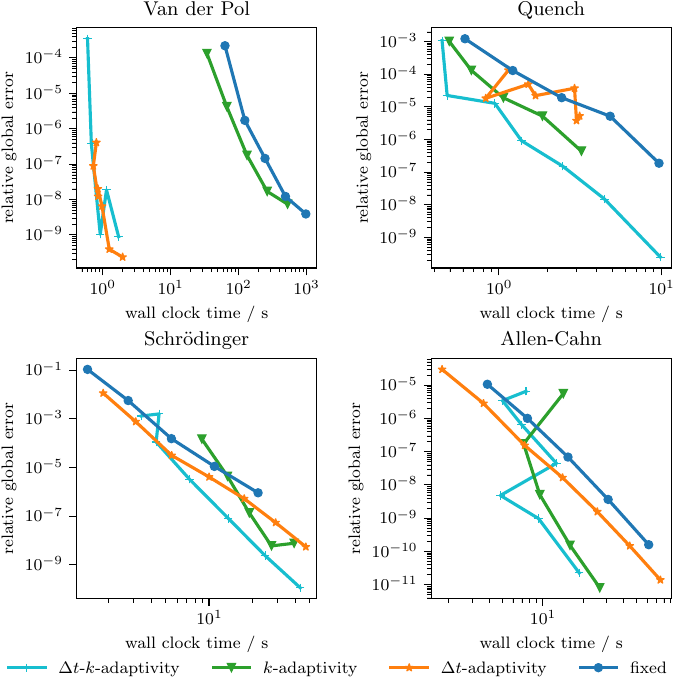}
    \caption{\reviewcomment{Updated van der Pol panel with $\mu=1000$ and Allen-Cahn panel with forcing term.}Wall-clock time versus the global error relative to the magnitude of the solution. For $\Delta t$- and $\Delta t$-$k$-adaptivity, the tolerance for the local error is adjusted, whereas the step size controls the accuracy of the other strategies. Note that for the Quench problem in $\Delta t$-adaptivity, the global error is not necessarily reduced when choosing a smaller tolerance because the error depends sensitively on the way the transition to runaway heating is resolved.\deleted{As the time-scale changes only little in Allen-Cahn, we can see that the fixed- and the $\Delta t$-adaptivity algorithm perform virtually identically with no meaningful overhead due to adaptivity.} \added{Adaptive resolution vastly enhances efficiency in all cases.}}
    \label{fig:time-error}
\end{figure}
\autoref{fig:time-error} shows error versus wall-clock time for the four problems from~\autoref{sec:Problems}.
\added{Choosing the step size adaptively is beneficial in all cases, but it proves particularly essential for the van der Pol case.
$\Delta t$-$k$-adaptive SDC is particularly efficient for the PDE examples.
However, the IMEX-scheme used in Allen-Cahn becomes unstable for large step sizes, causing around 100 restarts for loose tolerances and decreasing the efficiency of $\Delta t$-$k$-adaptive SDC.
A remedy is to limit the step size, which proved to increase efficiency in these cases in tests not shown here.
Note that even though efficiency is decreased, $\Delta t$-$k$-adaptivity provides physical solutions for any tolerance here.
}
\deleted{For all three PDE examples, Quench, Schr{\"o}dinger and Allen-Cahn, $\Delta t$-$k$-adaptive SDC is the most efficient variant, although $\Delta t$-adaptive SDC is very close for the Schr{\"o}dinger equation.
Only for the mildly stiff van der Pol ODE example is $\Delta t$-adaptive SDC is slightly more efficient.
This confirms that choosing both time step and number of iterations adaptively in SDC is generally a good approach.}

\begin{figure}
    \centering
    \includegraphics{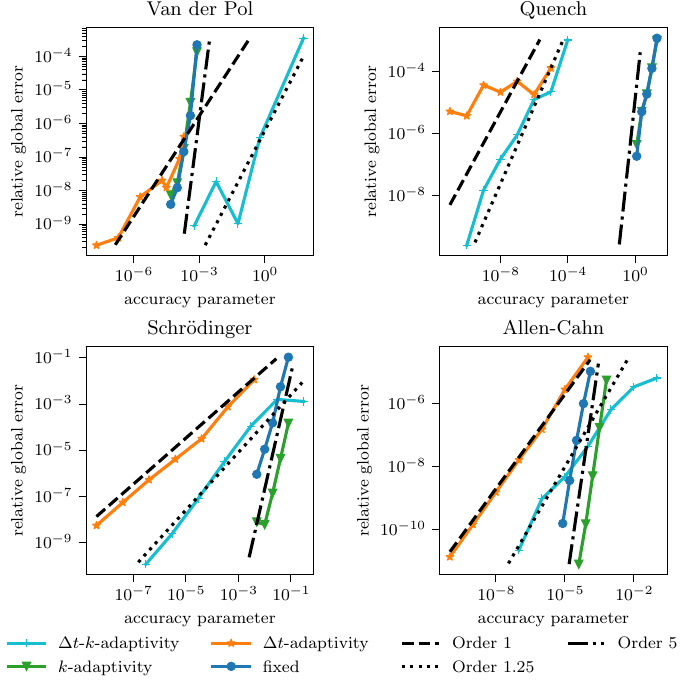}
    \caption{\reviewcomment{This figure is new}\added{Accuracy versus the parameter used to control it. For the fixed and $k$-adaptivity schemes this means step size and  $\epsilon_\mathrm{TOL}$ for step size-adaptive variants. The methods all have a global order of accuracy of 5 with respect to step size. We expect $e\propto \epsilon_\mathrm{TOL}$ for $\Delta t$-adaptivity and $e\propto \epsilon_\mathrm{TOL}^{1.25}$ for $\Delta t$-$k$-adaptivity, which we also observe in practice. Only in the combination of $\Delta t$-adaptivity and Quench, do we see significant deviations due to the complicating effect of the transition to runaway heating.}}
    \label{fig:param-error}
\end{figure}

\added{We investigate the relationship between $\epsilon_\mathrm{TOL}$ and the resulting error in \autoref{fig:param-error}.}
\added{We select the step size by controlling the local error of a lower order method, see \autoref{eq:error_estimate} and \autoref{eq:step_size_update}.
The resulting scaling for the global error is $e\propto \epsilon_\mathrm{TOL}^{q/p+1}$, with $p$ and $q$ as in \autoref{eq:error_estimate}.
That means we expect a linear dependence of global error on step size for $\Delta t$-adaptivity for any $q$.
For $\Delta t$-$k$-adaptivity, on the other hand, the scaling depends on the choice of quadrature; in our case the result is $e\propto \epsilon_\mathrm{TOL}^{5/4}$.
We observe this scaling in practice to a reasonable degree.
Some deviations due to stiffness or stability restrictions are not unexpected and occur also in embedded RKM.
}

To mitigate the cost of restarts, we use interpolation of the polynomial to obtain an initial guess after the collocation problem has converged with too large a step size in $\Delta t$-$k$-adaptivity, as discussed in \autoref{subsec:adaptive_time-stepping}. This indeed reduces the number of iterations in our numerical tests, but in the case of diagonally preconditioned SDC requires all-to-all communication, incurring additional cost.
However, the number of such restarts is small for all problems under consideration, such that the overall computational cost changes only little.
For problems with dynamics that lead to more restarts, on the other hand, interpolation can reduce the computational cost more significantly.

\begin{figure}
    \centering
    \includegraphics{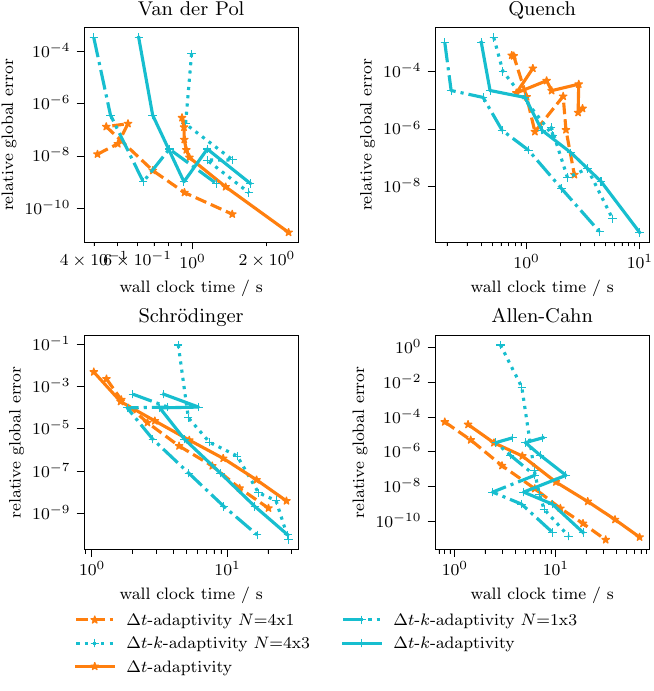}
    \caption{\reviewcomment{Updated van der Pol panel with $\mu=1000$ and Allen-Cahn panel with forcing term}Timings for parallel flavors of step size adaptive SDC. The number of processes is denoted, if larger than one, as $N$=[number of processes for GSSDC]x[number of processes for diagonal SDC]. For van der Pol and Allen-Cahn, GSSDC gives decent parallel efficiency, whereas better parallel efficiency is observed with diagonal SDC for the other problems. \replace{For the PDEs, w}{W}all time is roughly halved when using three processes in diagonal SDC. \added{Stability issues in the IMEX solver of Allen-Cahn can lead to very large errors or crashes when combining diagonal SDC and GSSDC.}}
    \label{fig:parallelspeedup}
\end{figure}

\paragraph{Parallel variants of SDC}
\autoref{fig:parallelspeedup} shows error against wall clock time for two serial adaptive variants of SDC, a $\Delta t$-adaptive block parallel variant of SDC, a $\Delta t$-$k$-adaptive variant with diagonal preconditioner and a $\Delta t$-$k$-adaptive combination of the latter two.

For the three PDE examples, the $\Delta t$-$k$-adaptive SDC with diagonal preconditioner ("parallel-across-the-method")~\cite{GT_precons} is the most efficient variant.
For the ODE \added{and Allen-Cahn} example\added{s}, the $\Delta t$-adaptive block parallel GSSDC variant \replace{performs best}{also shows good speedup}.
\deleted{The reason is that in the PDE case, right-hand side evaluations and implicit solves make up a larger fraction of the compute time, and these are performed in parallel in diagonal SDC.}
\replace{As Schr{\"o}dinger and Allen-Cahn are solved with spectral methods, the implicit solves are relatively cheap and the parallel efficiency is not as good as for the Quench problem.}{Since the implicit solvers are relatively expensive in the Quench problem, which uses a Newton scheme and finite difference discretization, parallel efficiency of diagonal SDC is especially good.}
Note that the serial comparison runs \added{for diagonal SDC} were also performed with diagonal preconditioners as these were found to give the best performance in our tests.

Unfortunately, combining both diagonal SDC and GSSDC does not further improve performance.
The reason seems to be that the performance of GSSDC is highly sensitive to the preconditioner and that it works best with implicit Euler for reasons not yet well understood.
Switching to a diagonal preconditioner in GSSDC significantly increases the number of iterations required which negates any performance gains from parallelization.
A deeper investigation of this issue is left for future work.

\subsection{Comparing efficiency against Runge-Kutta methods}
\label{sec:RKcomp}
\begin{figure}
    \centering
    \includegraphics{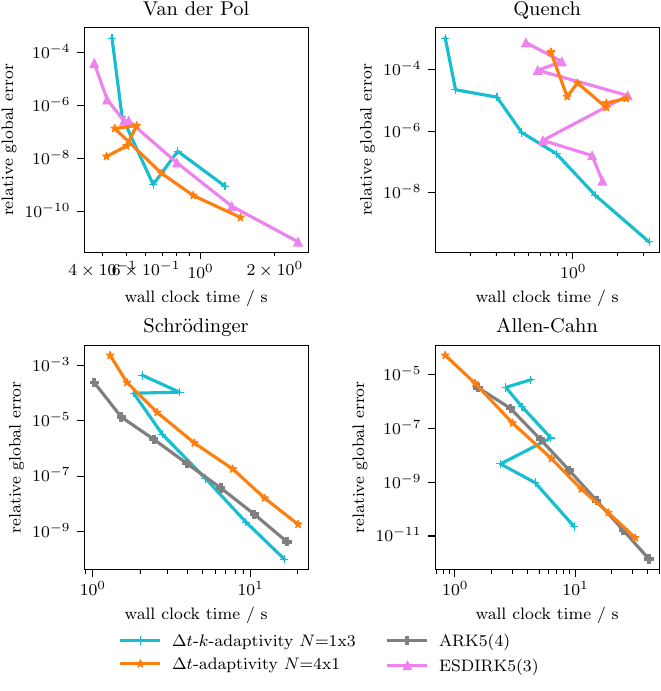}
    \caption{\reviewcomment{Updated van der Pol panel with $\mu=1000$ and Allen-Cahn panel with forcing term}Comparison of embedded RKM with parallel versions of $\Delta t$- and $\Delta t$-$k$-adaptive SDC from \autoref{fig:parallelspeedup}. All methods are of order 5. \deleted{As the van der Pol example is not too stiff, the explicit RKM performs best. For the very stiff PDEs, on the other hand, SDC can outperform the RKM. We do not show the explicit RKM for these problems because it ran into stability issues, producing unphysical solutions or no solutions at all within a reasonable time frame relative to the implicit methods. } \added{We find that we can outperform the RKM with one version of step size adaptive and parallel SDC for all problems.}\deleted{ Note that the ARK5(4) method performs very similar for Allen-Cahn as serial $\Delta t$-$k$-adaptivity, but the parallel version is much faster.}}
    \label{fig:RK-comparison}
\end{figure}

After establishing that adaptive SDC is generally more efficient than SDC with fixed $\Delta t$ and $k$, we now compare the run-times of $\Delta t$- and $\Delta t$-$k$-adaptive SDC against state-of-the-art embedded RKM of the same order.
\autoref{fig:RK-comparison} shows error against wall clock time  for two variants of SDC against different embedded RKM for our four benchmark problems.
\deleted{For the van der Pol problem, we compare against Cash-Karp's method  cite cash\_carp , an explicit pair of orders 4 and 5.}
For \added{van der Pol and }Quench, we compare against ESDIRK5(3)~\cite{ESDIRK}, a singly diagonally implicit, stiffly accurate pair of orders 3 and 5 with an explicit first stage.
For Schr{\"o}dinger and Allen-Cahn, we compare against ARK5(4)~\cite{KENNEDY2019183}, an additive pair of a singly diagonally implicit stiffly accurate L-stable embedded method of orders 4 and 5 and an explicit embedded method of orders 4 and 5.
\added{All problems are sufficiently stiff such that explicit RKM were encountering stability issues, produced unphysical solutions or no solutions at all within a reasonable time frame relative to the implicit methods.}

\deleted{Since the van der Pol problem is only mildly stiff, SDC cannot beat the explicit CP5(4) RKM.
For the three PDE examples, however, $\Delta t$-$k$-adaptive SDC with parallel preconditioner is always at least competitive, providing similar or better performance than RKM.}
\added{We find that at least one type of step size adaptive parallel SDC can outperform the RKM for all examples.
In our testing, the RKM is only faster for Schr{\"o}dinger when low accuracy is sufficient.}
Therefore, adaptive SDC is not only more efficient than standard SDC but also competitive with state-of-the-art adaptive RKM\deleted{, at least for the PDE case}.
SDC also offers easy tuning of the order, simply by changing parameters.
By contrast, ARK5(4) is the highest order additive RKM we are aware of.
\added{In experiments not shown here, we repeated the comparison with third order accurate methods.
We found the same trends as in \autoref{fig:RK-comparison} with at least one SDC method outperforming the RKM, even though lower order SDC methods allow for less concurrency.}
\added{Note that there exist also stage-parallel RKM, which allow parallelism with a similar idea as diagonal SDC, see e.g.~\cite{PaRK}.
However, they often have worse stability than their serial counterparts and we were unable to find an embedded method of order higher than 4 in the literature.}

\reviewcomment{deleted numerical results on resilience section}

\section{Summary}\label{sec:Discussion}
The paper adopts concepts from adaptive embedded Runge-Kutta methods to spectral deferred correction\deleted{s}.
We propose procedures to control the step size $\Delta t$ in SDC, the iteration number $k$ or both.
Our adaptive techniques can also be used for variants of SDC that are parallelizable, either in a Gauss-Seidel "parallel-across-the-steps" fashion (GSSDC) or in a "parallel-across-the-method" way by using a diagonal preconditioner.
Numerical examples demonstrate that adaptive SDC is more efficient than SDC with fixed step size and iteration number and that adaptive parallel SDC can be competitive with embedded Runge-Kutta methods for the integration of \replace{three}{four} complex, nonlinear time-dependent \replace{PDEs}{problems}.

One advantage of SDC is flexibility due \added{to }its iterative nature, making it easy to include techniques like splitting.
The $\Delta t$-$k$-adaptive algorithm, in particular, retains much of the flexibility of SDC.
Preconditioners or reduced accuracy spatial solves tailored to specific problems can be used.
Also, spatial adaptivity may be leveraged in unique ways~\cite{space_adaptive_SDC_cardiac, WeiserChegini2022, WeiserCheginiPreprint}.

Adaptive time stepping has been explored for revisionist integral deferred correction\deleted{s} (RIDC)~\cite{christlieb_revisionist_2015}, which is similar to GSSDC.
However, they find that the increment is a poor choice for an error estimate, because of the accumulation of local errors from different steps in the increment.
The reason why it works well for GSSDC is that instead of allowing maximal pipelining, GSSDC solves only fixed size blocks of steps at a time, each with the same step size.
While sacrificing flexibility compared to RIDC, this allows us to view the accumulated errors inside the blocks as a local error of the block and enables efficient step size selection in GSSDC.
More elaborate step size selection for GSSDC that allow different step sizes for each step in the blocks like in RIDC are left for future work.

\added{The similarity of our methods to embedded RKM means any step size controller used in that context can also be used in SDC.
It is well known that ideas like limiting the step size can further boost efficiency on a problem specific basis.
So while our method only uses a single tolerance parameter, more elaborate step size update equations can readily be employed by domain scientists with intimate knowledge of the problem at hand.}

\deleted{We also demonstrate that adaptive SDC provides algorithm-based resilience against bitflips in the solution, similar to what has been documented for adaptive RKM.
Fault tolerance based on algorithms may yield significant advantages as supercomputers continue to grow in complexity.
While adaptive SDC can recover a large range of faults, we found that we are unable to recover from faults in the initial condition of a time step.
These could be further protected on the application level, for instance by replication, or they could be stored in a region of memory protected by error correction codes that can detect and recover from faults, whereas all other stages of the collocation problem may be stored in unprotected regions of memory.
This would use available memory more efficiently, in particular in terms of bandwidth, by avoiding the checksums that are involved in hardware error correction codes.
Exploration of such more sophisticated resilience schemes is left for future work.}

\paragraph*{Acknowledgements}
The authors gratefully acknowledge computing time granted on JUSUF through the CSTMA project at J{\"u}lich Supercomputing Centre.
\added{The authors also gratefully acknowledge the helpful comments from the reviewers.}

\paragraph*{Funding}
We thankfully acknowledge funding from the European High-Performance Computing Joint Undertaking (JU) under grant agreement No 955701. 
The JU receives support from the European Union’s Horizon 2020 research and innovation programme and Belgium, France, Germany, and Switzerland. 
We also thankfully acknowledge funding from the German Federal Ministry of Education and Research (BMBF) grant 16HPC048.

\paragraph*{Availability of supporting data}
For instructions on how to reproduce the results shown in this paper, please consult the project section in the documentation of the pySDC code at~\url{https://parallel-in-time.org/pySDC/projects/Resilience.html}.

\paragraph*{Ethical approval}
Not applicable

\paragraph*{Competing interests}
All authors certify that they have no affiliations with or involvement in any organization or entity with any financial interest or non-financial interest in the subject matter or materials discussed in this manuscript.

\paragraph*{Author contributions}
T.S. wrote the main manuscript text and carried out the implementations.
T.L. prepared Figure 1.
All authors reviewed the manuscript and considerably contributed to both text and results.

\bibliography{bibliography.bib}


\end{document}